\documentclass[twoside,a4paper,12pt,centertags,reqno]{amsart} 
\usepackage{amsmath,amssymb,verbatim,vmargin}
\usepackage{color}
\usepackage{tikz}
\usepackage{tikz-cd}
\usetikzlibrary{matrix,calc}
\usepackage{hyperref}
\hypersetup{colorlinks,bookmarks=true,linktocpage=true,citecolor=blue,linkcolor=magenta}

\usepackage{eulervm}   

\allowdisplaybreaks 
\usepackage{color}
\usepackage[all]{xy} 

\usepackage{mathtools}
\usepackage{MnSymbol} 

\theoremstyle{plain}
\newtheorem{thm}{Theorem}[section]
\newtheorem{lem}[thm]{Lemma}

\theoremstyle{definition}

\newtheorem{rem}[thm]{Remark}

\usepackage{enumerate}
\usepackage{enumitem}

\newcommand{\bR}{{\mathbb R}}

\newcommand{\cH}{{\mathcal H}}

\newcommand{\fD}{{\mathbf D}}

\def\barint_#1{\mathchoice
            {\mathop{\vrule width 6pt
height 3 pt depth -2.5pt
                    \kern -9.5pt
\intop \kern -4pt}\nolimits_{#1}}%
            {\mathop{\vrule width 5pt height
3 pt depth -2.6pt
                    \kern -6.5pt
\intop \kern -4pt}\nolimits_{#1}}%
            {\mathop{\vrule width 5pt height
3 pt depth -2.6pt
                    \kern -6pt
\intop \kern -4pt}\nolimits_{#1}}%
            {\mathop{\vrule width 5pt height
3 pt depth -2.6pt
          \kern -6pt \intop \kern -4pt}\nolimits_{#1}}}
          
           \def\bariint_#1{\mathchoice
            {\mathop{\vrule width 15pt
height 3 pt depth -2.5pt
                    \kern -15.8pt
\intop \kern -8pt\intop \kern -4pt}\nolimits_{#1}}%
            {\mathop{\vrule width 9pt height
3 pt depth -2.6pt
                    \kern -10.5pt
\intop \kern -8pt\intop \kern -4pt}\nolimits_{#1}}%
            {\mathop{\vrule width 9pt height
3 pt depth -2.6pt
                    \kern -10pt
\intop \kern -8pt\intop \kern -4pt}\nolimits_{#1}}%
            {\mathop{\vrule width 9pt height
3 pt depth -2.6pt
          \kern -8pt \intop \kern -10pt\intop \kern -4pt}
      \nolimits_{  #1}}}

\def\barintlim_#1{\mathchoice
            {\mathop{\vrule width 6pt
height 3 pt depth -2.5pt
                    \kern -8.8pt
\intop \kern -4pt}\limits_{#1}}%
            {\mathop{\vrule width 5pt height
3 pt depth -2.6pt
                    \kern -6.5pt
\intop \kern -4pt}\limits_{#1}}%
            {\mathop{\vrule width 5pt height
3 pt depth -2.6pt
                    \kern -6pt
\intop \kern -4pt}\limits_{#1}}%
            {\mathop{\vrule width 5pt height
3 pt depth -2.6pt
          \kern -6pt \intop \kern -4pt}\limits_{#1}}}
          
           \def\bariintlim_#1{\mathchoice
            {\mathop{\vrule width 15pt
height 3 pt depth -2.5pt
                    \kern -15.8pt
\intop \kern -8pt\intop \kern -4pt}\limits_{#1}}%
            {\mathop{\vrule width 9pt height
3 pt depth -2.6pt
                    \kern -10.5pt
\intop \kern -8pt\intop \kern -4pt}\limits_{#1}}%
            {\mathop{\vrule width 9pt height
3 pt depth -2.6pt
                    \kern -10pt
\intop \kern -8pt\intop \kern -4pt}\limits_{#1}}%
            {\mathop{\vrule width 9pt height
3 pt depth -2.6pt
          \kern -8pt \intop \kern -10pt\intop \kern -4pt}
      \limits_{  #1}}}
          
\renewcommand{\iint}{\int \kern -8pt\int}       

\newcommand{\RE}{\text{{\rm Re}}\,}

\numberwithin{equation}{section}
\makeatletter
\@namedef{subjclassname@2020}{\textup{2020} Mathematics Subject Classification}
\makeatother

\title{On Landau-Kato inequalities via semigroup orbits}

\author{Yi C. Huang} 
\address{School of Mathematical Sciences, Nanjing Normal University, Nanjing 210023, People's Republic of China}
\email{Yi.Huang.Analysis@gmail.com}
\urladdr{https://orcid.org/0000-0002-1297-7674}

\author{Yanlu Lian} 
\address{School of Mathematics, Hangzhou Normal University, Hangzhou 311121, People's Republic of China}
\email{yllian@hznu.edu.cn}

\author{Fei Xue} 
\address{School of Mathematical Sciences, Nanjing Normal University, Nanjing 210023, People's Republic of China}
\email{05429@njnu.edu.cn}

\date{\today} 

\thanks{Research of the authors is supported by the National NSF grants of China (nos. 11801274 and 12201307), the Jiangsu Provincial NSF grant (no. BK20210555), 
the China Scholarship Council, and the Scientific Research Fund of Zhejiang Provincial Education Department (no. Y202250092).
YCH thanks Gerd Herzog (KIT) and Masayuki Hayashi (Univ. of Pisa) for helpful communications.}

\begin{document}

\begin{abstract}
Let $\omega>0$.
Given a strongly continuous semigroup $\{e^{tA}\}$ on a Banach space and an element $f\in\fD(A^2)$ satisfying the exponential orbital estimates
$$\|e^{tA}f\|\leq e^{-\omega t}\|f\| \quad\text{and}\quad \|e^{tA}A^2f\|\leq e^{-\omega t}\|A^2f\|,\quad t\geq0,$$
a dynamical inequality for $\|Af\|$ in terms of $\|f\|$ and $\|A^2f\|$ was derived by Herzog and Kunstmann (Studia Math., 2014).
Here we provide an improvement of their result by relaxing the exponential decay to quadratic,
together with a simple and direct way recovering the usual Landau inequality.
Herzog and Kunstmann also demanded an analogue, again via semigroup orbits, for the Kato type inequality on Hilbert spaces.
We provide such a result by using Hayashi-Ozawa machinery [Proc. Amer. Math. Soc., (2017)] which in turn relies on Hilbertian geometry.
\end{abstract}

\maketitle


\section{Introduction}

The classical Landau inequality for functions defined on $\bR_+:=(0,\infty)$ states
$$\|f'\|_\infty^2\leq4\|f\|_\infty\|f''\|_\infty,$$
where $f'=\frac{df}{dx}$ and $\|\cdot\|_\infty=\|\cdot\|_{L^\infty(\bR_+)}$.
The constant 4 can be replaced by 2 if one considers the Hilbert space $L^2(\bR_+)$ instead of $L^\infty(\bR_+)$; 
this is a classical result of Hardy, Littlewood, and P\'olya.
Let $A$ be the generator of a strongly continuous contractive\footnotemark
\footnotetext{Namely, $\|e^{tA}\|_{X\rightarrow X}\leq1$ for all $t\geq0$.} semigroup on a Banach space $X$, with norm $\|\cdot\|=\|\cdot\|_X$ for notational convenience.
If $f\in\fD(A^2)\backslash\{0\}$, we also have the abstract Landau inequality
\begin{equation} \label{eqn:localLandau}
\|Af\|^2\leq4\|f\|\|A^2f\|.
\end{equation}
If $X$ is in addition a Hilbert space, the abstract Hardy-Littlewood-P\'olya inequality
\begin{equation} \label{eqn:Kato}
\|Af\|^2\leq2\|f\|\|A^2f\|
\end{equation}
was obtained by Kato in \cite{Kat71} (see also Hayashi and Ozawa \cite{HayOza17} for an elegant proof based on Hilbertian geometry).
For related results and a recent survey of such norm inequalities, 
see for example Kallman-Rota \cite{KalRot70}, Gindler-Goldstein \cite{GinGol75}, Chernoff \cite{Che79}, Kwong-Zettl \cite{KwoZet92}, 
and Gesztesy-Nichols-Stanfill \cite{GezNicSta21}.

Let $0<\phi\in C^2(\bR_+)$ be a decreasing convex function so that 
 $$\phi(0)=-\phi'(0)=1\quad\text{and}\quad \phi(t)\leq1-t+\frac{t^2}{2}.$$
Typical candidates for $\phi(t)$ are 
$$e^{-t} \qquad \text{and}\qquad\frac{1}{1+t+\frac{t^2}{2}}.$$
Assuming proper evolutionary bounds (measured by above function $\phi$) along two ``semigroup orbits",
the following theorem extends and improves the corresponding result of Herzog and Kunstmann in \cite{HerKun14} for rigid exponential decay.

\begin{thm} \label{thm:localLandau}
Let $\{e^{tA}\}_{t\geq0}$ be a strongly continuous semigroup on a Banach space $X$.
Suppose that for some $f\in\fD(A^2)\backslash\{0\}$ and $\omega>0$, and for all $t\geq0$,
\begin{equation} \label{eqn:orbital}
\|e^{tA}f\|\leq \phi(\omega t)\|f\|\quad\text{\,\, and\,\,}\quad \|e^{tA}A^2f\|\leq \phi''(\omega t)\|A^2f\|.
\end{equation}
Then, for the quantities 
\begin{equation} \label{eqn:abc}
a=\|f\|,\quad b=\frac{\|Af\|}{\omega}\quad\text{and}\quad c=\frac{\|A^2f\|}{\omega^2},
\end{equation}
we have
\begin{equation} \label{eqn:expineq}
\phi(s)(a+c)+s(c-b)+(a-c)\geq0,\quad\forall\,s\geq0.
\end{equation}
In particular, 
\begin{equation} \label{eqn:localLandau'}
b^2\leq 4ac.
\end{equation}
\end{thm}

\begin{rem}
When $\phi(t)=e^{-t}$, the dynamical inequality \eqref{eqn:expineq} was obtained by Herzog and Kunstmann in \cite{HerKun14}.
Moreover, the following relation was justified:
\begin{equation} \label{eqn:order}
a\leq b\leq c.
\end{equation}
In view of the bounds \eqref{eqn:orbital}, it is natural to require that $\phi$ is decreasing and convex.
\end{rem}

\begin{rem}
The precision of polynomial decay for operator semigroups is important in applications to evolution equations,
see e.g. Borichev-Tomilov \cite{BorTom10}.
\end{rem}

Note that $b$ and $c$ depend on the locality parameter $\omega$,
although it is factored out in $b^2\leq 4ac$.
Thus \eqref{eqn:expineq}-\eqref{eqn:localLandau'} can be understood as the $\omega$-localised version of \eqref{eqn:localLandau}.

\section{Proof of \eqref{eqn:expineq}}

We follow the beautiful proof of \eqref{eqn:expineq} given by Herzog and Kunstmann in \cite{HerKun14} when $\phi(t)=e^{-t}$.
For $f\in\fD(A^2)\backslash\{0\}$ we have the following integral representation 
$$\begin{aligned}
e^{tA}f&=f+\int_0^te^{\tau A}Afd\tau\\
&=f+\int_0^t\left(Af+\int_0^\tau e^{\tau' A}A^2fd\tau'\right)d\tau\\
&=f+tAf+\int_0^t(t-\tau)e^{\tau A}A^2fd\tau,\quad t\geq0.
\end{aligned}$$
In the last step we used the integration by parts.
Then by the triangle inequality and the orbital decay assumptions for $e^{tA}$ on $f$ and $A^2f$, we have
\begin{equation} \label{eqn:triangle}
\begin{aligned}
\phi(\omega t)\|f\|&\geq \|e^{tA}f\|\\
&\geq\|f+tAf\|-\left(\int_0^t(t-\tau)\phi''(\omega \tau)d\tau\right)\|A^2f\|\\
&\geq t\|Af\|-\|f\|-\|A^2f\|\left(\frac{\phi(\omega t)-1}{\omega^2}+\frac{t}{\omega}\right).
\end{aligned}
\end{equation}
In the last step we used again the integration by parts, together with $\phi(0)=-\phi'(0)=1$.
Rewriting \eqref{eqn:triangle} using the notations \eqref{eqn:abc} then gives \eqref{eqn:expineq}.

\section{A simple and direct proof of \eqref{eqn:localLandau'} via \eqref{eqn:expineq}}

By homogeneity in \eqref{eqn:expineq} and \eqref{eqn:localLandau'}, we can take $a=1$.
Re-organising \eqref{eqn:expineq} and using the assumptions on $\phi$ (namely, $\phi\leq\phi(0)=1$ and $\phi(t)\leq1-t+\frac{t^2}{2}$), we have
$$\begin{aligned}
b&\leq\frac{(\phi(s)+1)+(\phi(s)-1+s)c}{s}\\
&\leq\frac2s+\frac{sc}{2},\quad\forall\,s\geq0.
\end{aligned}$$
Minimising the right hand side we get $b\leq2\sqrt c$, which is \eqref{eqn:localLandau'} for $a=1$.

\begin{rem}
In above reasoning, the order relation \eqref{eqn:order} is not needed.
\end{rem}

\begin{rem}
In \cite{HerKun14} with $\phi(t)=e^{-t}$, an intermediate (or interpolated) estimate 
\begin{equation} \label{eqn:localLandau''}
b\leq c-(a+c)g^{-1}\left(\frac{c-a}{c+a}\right),
\end{equation}
where $$g(x)=\begin{cases}x(1-\log x),& x\in(0,1],\\ 0, &x=0,\end{cases}$$
was obtained via minimising the left hand side of \eqref{eqn:expineq}.
However, to evaluate the quality of \eqref{eqn:localLandau''}, 
a rather involved limiting argument was used to recover \eqref{eqn:localLandau'}.
\end{rem}

\section{Kato inequality via semigroups orbits}

As Open Problem, Herzog and Kunstmann demanded at the end of \cite{HerKun14} an analogue of Theorem \ref{thm:localLandau} for \eqref{eqn:Kato}.
The aim of this section is to provide such a result by using the Hilbertian machinery of Hayashi and Ozawa in \cite{HayOza17}.

Thus we are working in a Hilbert space $\cH$ with norm $\|\cdot\|$ and inner product $(\cdot\mid\cdot)$.
The following lemma based on Hilbertian geometry is elementary but very useful. 

\begin{lem} \label{lem:Hilbert}
For all $f_1$, $f_2$ and $f_3$ that are elements in $\cH$, we have
$$\|f_1+f_2+f_3\|^2+\|f_2\|^2=\|f_1\|^2+\|f_3\|^2+2\RE(f_1+f_2\mid f_2+f_3).$$
\end{lem}

\begin{proof}
By straightforward expansion of the involved norms and inner products.
\end{proof}

For any $\lambda\geq0$ and $f\in\fD(A^2)\backslash\{0\}$, taking in Lemma \ref{lem:Hilbert}
$$f_1=A^2f, \quad f_2=\lambda Af\quad \text{and}\quad f_3=\lambda^2f$$
and dropping out the nonnegative term $\|f_1+f_2+f_3\|^2$, we have\footnotemark
\footnotetext{Note that the case $\lambda=0$ for \eqref{eqn:Katoc} is trivial.} 
\begin{equation} \label{eqn:Katoc}
\lambda^2\|Af\|^2\leq\|A^2f\|^2+\lambda^4\|f\|^2
\end{equation}
once we assume the following ``restricted dissipativity"
\begin{equation} \label{eqn:orbidiss}
\RE(A(A+\lambda)f\mid(A+\lambda)f)\leq0,\quad \lambda>0.
\end{equation}
Note that the inequality \eqref{eqn:Katoc} is equivalent to the Kato inequality \eqref{eqn:Kato}.
 
Now, suppose for some $f\in\fD(A^2)\backslash\{0\}$ one has a family of orbital estimates
\begin{equation} \label{eqn:orbitalc1}
\|e^{tA}(A+\lambda)f\|\leq \|(A+\lambda)f\|, \quad t\geq0, \quad \lambda>0.
\end{equation}
From \eqref{eqn:orbitalc1} we derive
$$\RE (e^{tA}(A+\lambda)f\mid (A+\lambda)f)\leq ((A+\lambda)f\mid (A+\lambda)f)$$
and its infinitesimal version \eqref{eqn:orbidiss}. The Kato inequality \eqref{eqn:Kato} then follows.

Hence we proved the following ``semigroup orbits" result in the spirit of \cite{HerKun14}.

\begin{thm} \label{thm:localKato}
Let $\{e^{tA}\}_{t\geq0}$ be a strongly continuous semigroup on a Hilbert space $\cH$ with norm $\|\cdot\|$.
Then, the family of semigroup orbital estimates \eqref{eqn:orbitalc1} for some element $f\in\fD(A^2)\backslash\{0\}$
implies the Kato inequality \eqref{eqn:Kato} for the same $f$.
\end{thm}

\begin{rem}
The orbital estimates \eqref{eqn:orbital} require an exponential decay of the semigroup bound on $A^2f$ and $f$ (and by homogeneity, on $\lambda f$ for all $\lambda>0$),
whereas \eqref{eqn:orbitalc1} requires only the contractivity of the semigroup on $(A+\lambda) f$ for all $\lambda>0$.
\end{rem}

\bigskip

\section*{\textbf{Compliance with Ethical Standards}}

\bigskip

\textbf{Conflict of interest} The authors have no known competing financial interests 
or personal relationships that could have appeared to influence this reported work.

\bigskip

\textbf{Availability of data and material} Not applicable.

\bigskip

\bibliographystyle{alpha}
 
\bibliography{HuaY-LiaY-XueF-LanKatSemigroupOrbits}

\end{document}